\documentclass[12pt]{article}

\usepackage{amsmath}
\usepackage{amssymb}
\usepackage{dsfont} 

\usepackage{graphicx}

\usepackage{color}

\newcommand{\be}{\begin{equation}}
\newcommand{\ee}{\end{equation}}
\newcommand{\bel}[1]{\begin{equation}\label{#1}}
\newcommand{\bea}{\begin{eqnarray}}
\newcommand{\eea}{\end{eqnarray}}
\newcommand{\balign}{\begin{align}}
\newcommand{\ealign}{\end{align}}
\newcommand{\ba}{\begin{array}}
\newcommand{\ea}{\end{array}}
\newcommand{\bfig}{\begin{figure}}
\newcommand{\efig}{\end{figure}}

\newcommand{\eref}[1]{(\ref{#1})}

\allowdisplaybreaks

\newcommand{\bra}[1]{\mbox{$\langle \, {#1}\, |$}}
\newcommand{\ket}[1]{\mbox{$| \, {#1}\, \rangle$}}
\newcommand{\exval}[1]{\mbox{$\langle \, {#1}\, \rangle$}}
\newcommand{\inprod}[2]{\mbox{$\langle \, {#1} \, | \, {#2} \, \rangle$}}
\newcommand{\Prob}[1]{\mbox{${\rm Prob}\left[ \, {#1}\, \right]$}}

\newcommand{\bfx}{\mathbf{x}}
\newcommand{\bfy}{\mathbf{y}}
\newcommand{\bfz}{\mathbf{z}}

\newcommand{\bfr}{\mathbf{r}}
\newcommand{\bfs}{\mathbf{s}}

\newcommand{\rme}{\mathrm{e}}

\newcommand{\half}{\frac{1}{2}}

\newcommand{\comm}[2]{\mbox{$[\,{#1}\,,\,{#2}\,]$}}

\newcommand{\C}{{\mathbb C}}
\newcommand{\R}{{\mathbb R}}

\newcommand{\Z}{{\mathbb Z}}

\renewcommand{\S}{\mathbb S}

\newcommand{\bfeta}{\boldsymbol{\eta}}

\newtheorem{theo}{Theorem}[section]
\newtheorem{lmm}[theo]{Lemma}
\newtheorem{df}[theo]{Definition}

\newtheorem{cor}[theo]{Corollary}
\newtheorem{rem}[theo]{Remark}
\def\qed{\hfill$\Box$\par\medskip\par\relax}

\begin{document}

\title{Self-Duality for the Two-Component Asymmetric Simple Exclusion Process}
\author{V. Belitsky$^{1}$ 
\and G.M.~Sch\"utz$^{2,3}$ 
}

\maketitle

{\small
\noindent $^{~1}$Instituto de Matem\'atica e Est\'atistica,
Universidade de S\~ao Paulo, Rua do Mat\~ao, 1010, CEP 05508-090,
S\~ao Paulo - SP, Brazil
\\
\noindent Emails: belitsky@ime.usp.br  

\smallskip
\noindent $^{~2}$Institute of Complex Systems II,
Forschungszentrum J\"ulich, 52425 J\"ulich, Germany
\\
\noindent Email: g.schuetz@fz-juelich.de

\smallskip
\noindent $^{~3}$Interdisziplin\"ares Zentrum f\"ur Komplexe Systeme, Universit\"at
Bonn, Br\"uhler Str. 7, 53119 Bonn, Germany\\
\noindent URL: http://www.izks.uni-bonn.de
}

\begin{abstract}
We study a two-component asymmetric simple exclusion process (ASEP) that is equivalent to the
ASEP with second-class particles. We prove self-duality with respect to a family of 
duality functions which are shown to arise from the reversible measures of the process and the 
symmetry of the generator under the quantum algebra $U_q[\mathfrak{gl}_3]$.
We construct all invariant measures in explicit form and discuss some of their properties.
We also prove a sum rule for the duality functions.
\\[.3cm]\textbf{Keywords:} Asymmetric simple exclusion process; Second Class Particles; Duality;
Quantum algebras
\\[.3cm]\textbf{AMS 2000 subject classifications:} 82C20. Secondary: 60K35, 82C23
\end{abstract}

\newpage

\section{Introduction}

We consider an asymmetric simple exclusion process with two species of particles 
on the one-dimensional finite lattice $\Lambda = \{-L+1,\dots,L\}$. Its Markovian
dynamics can be described informally as follows.
Each site $i$ can be either empty (denoted by $0$) or occupied by at most one 
particle of type $A$ or of type $B$.
Thus we have local occupation numbers $\eta(k)\in\{A,0,B\}$.
We define the bonds $(k,k+1)$ of $\Lambda$ where $-L+1\leq k \leq L-1$.
Each bond carries a clock $i$ which rings independently of all other clocks 
after an exponentially distributed
random time with parameter $\tau_k$ where $\tau_k=r$ if $(\eta(k),\eta(k+1))\in 
\{(A,0),(0,B),(A,B)\}$ and $\tau_k=\ell$ if $(\eta(k),\eta(k+1))\in 
\{(0,A),(B,0),(B,A)\}$. When the clock rings the particle occupation variables are
interchanged and the clock acquires the corresponding new parameter.
Symbolically this process can be presented by the nearest neighbour particle
jumps
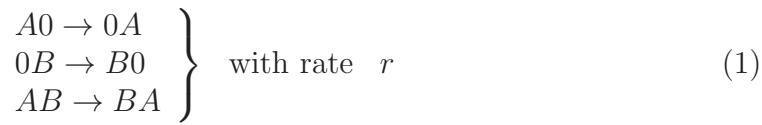
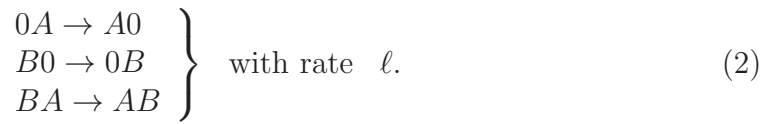
\bea
\label{2ASEPr}
\left. \ba{l}
A0 \to 0A \\
0B \to B0 \\
AB \to BA \ea \right\}
& \mbox{with rate} & r \\
\label{2ASEPl}
\left. \ba{l}
0A \to A0 \\
B0 \to 0B \\
BA \to AB \ea \right\}
& \mbox{with rate} & \ell .
\eea
We have reflecting boundary conditions, 
which means that no jumps from the left boundary site $-L+1$ to the left and no jumps from 
the right boundary site $L$ to the right are allowed.
We shall assume partially asymmetric hopping, i.e., $0 < r,\ell < \infty$. 
By interchanging the $B$-particles and vacancies this process turns into the
ASEP with second-class particles \cite{Ferr91}. 
We choose an even number of lattice sites exclusively for the
sake of convenience of notation. 

The objective of this work is to construct for the finite lattice in explicit form
all reversible measures and to prove self-duality with respect to a family of duality functions
that allows for the computation of expectations of the many-particle system in terms of
transition probabilities of the same process with only a small number of particles. 
It will transpire that this property, analogous to the well-known self-duality of the 
simple symmetric exclusion process
\cite{Ligg85}, arises from the fact proved in \cite{Beli15a} that the generator of this 
process commutes with a set of matrices which form a representation of the
quantum algebra
$U_q[\mathfrak{gl}(3)]$, which is the $q$-deformed universal enveloping algebra of 
the Lie algebra $\mathfrak{gl}(3)$ defined below (\eref{Uqglndef1} - \eref{UqglnSerre2}). 

The idea of deriving of duality relations from the representation matrices of
a non-abelian symmetry algebra of a 
generator of a Markov process goes back to Sch\"utz and 
Sandow \cite{Schu94} where this strategy
was applied to the symmetric partial exclusion process on arbitrary lattices. 
This is an interacting particle system with a $SU(2)$ symmetry where each lattice site 
can be occupied by
at most a finite number of particles. Next this symmetry approach was extended
to prove self-duality of the asymmetric simple exclusion process \cite{Schu97}, 
which is symmetric under the action of  quantum algebra
$U_q[\mathfrak{gl}(2)]$ and which is an integrable model solvable by Bethe ansatz. The self-duality
together with the integrability was used in \cite{Beli02} to study the time evolution 
of shock measures and in \cite{Imam11} to study current moments. By mapping the ASEP
to a lattice model of interface grwoth the duality function
can be interpreted as a lattice Cole-Hopf transformation \cite{Schu97} and is
therefore yields information on KPZ interface growth and the
moments of the partition function of a directed polymer \cite{LeDo12}.

The idea of using symmetries of the generator to obtain duality functions was employed 
again by Giardin\`a et. al. \cite{Giar07} to study heat conduction in the KMP model with 
$SU(1,1)$ symmetry and subsequently extended to other interacting particle systems, including
particle systems without conservation of particle number \cite{Giar09,Ohku10,Cari13,Cari15}.
Recently the $U_q[\mathfrak{gl}(2)]$ symmetry was extended to the non-integrable 
asymmetric generalization of the $SU(2)$-symmetric partial exclusion process \cite{Cari14}. 
Duality relations for new integrable models that can be solved
by Bethe ansatz and related methods were studied very recently in \cite{Boro14,Corw15}.

Here we prove self-duality for the two-component ASEP mentioned above
whose symmetry algebra $U_q[\mathfrak{gl}(3)]$ is larger than $SU(2)$, $SU(1,1)$ or their
$q$-deformations. We shall consider only finite systems, the construction and characterization
of the properties of the process on $\Z$ is out of the scope of this work.
The main novel feature is the presence of more than one conserved species
of particles. This leads to interesting non-local properties of the
duality functions and, through the integrability of the model, 
to the possibility of applications in the infinite volume limit 
employing exact computations along the lines of \cite{Chat10,Trac13}. 

The paper is structured as follows. In Sec. 2 we define the process and mention two results 
obtained in recent work \cite{Beli15a} that will be used here.
In Sec. 3 we state the main results of the present work. Sec. 4 is included for self-containedness. 
We describe some tools from linear algebra \cite{Lloy96,Schu01} used in the 
proofs, which are convenient, but not widely known in the probabilistic treatment of 
interacting particle systems.
In Sec. 5 we present the proofs of our results.

\section{The two-component ASEP}
\label{Sec:Definotat}

We define the process, introduce notation, and mention some results used in the proofs.

\subsection{State space and configurations}
\label{Sec:Defconfig}

It is convenient to introduce ternary local state variables $\eta(k) \in \mathbb{S}$ 
where $\mathbb{S} = \{0,1,2\}$.
We say that 0 represents occupation of a
site a particle of type $A$, 1 represents a vacant site and 2 represents occupation by a particle
of type $B$. Thus a configuration  is denoted by 
$\bfeta = \{ \eta(-L+1), \dots , \eta(L) \} \in \mathbb{S}^{2L}$.
We call this characterization of a configuration the occupation variable representation.
We shall repeatedly consider
configurations with a fixed number $N$ particles of type $A$ and $M$ particles of type $B$.
We denote configurations with this property by $\bfeta_{N,M}$ and the set of all such configurations
by $\mathbb{S}^{2L}_{N,M}$.

Equivalently we can specify a configuration $\bfeta$ uniquely by 
indicating the particle positions $\bfz$ on the lattice 
and write
\bel{ppeta}
\bfz = \{ \bfx, \bfy \}
\ee
with
\bel{pp}
\bfx := \{ x_i \, :\, \eta(x_i) = 0 \}, \quad \bfy := \{ y_i \, :\, \eta(y_i) = 2 \}
\ee
We call this the position representation.

Throughout this work we use the Kronecker-symbol defined by 
\be 
\delta_{\alpha,\beta} = \left\{ \ba{ll} 1 & \mbox{ if } \alpha=\beta \\ 0 & \mbox{ else } \ea \right.
\ee
for $\alpha,\beta$ from any set. We also introduce for $k,l \in \Lambda$
\bel{Thetaxz}
\Theta(k,l) := \left\{ \ba{ll} 1 \quad & k < l \\ 0 & k \geq l \ea \right. .
\ee
and the indicator function
\be
\mathbf{1}_A(\bfeta) = \left\{ \ba{ll} 1 & \mbox{ if } \bfeta \in A \\ 0 & \mbox{ else } \ea \right.
\ee
for subsets $A \subseteq \S^{2L}$. 
Some other functions of the configurations will play a role in our treatment:


\begin{df}
For $1\leq k < L$ we define the local permutation
\bel{etaexch}
\sigma^{kk+1}(\bfeta) = \{ \eta(-L+1), \dots \eta(k-1), \eta(k+1), \eta(k), \eta(k+2), \dots ,\eta(L) \}
=: \bfeta^{kk+1} .
\ee
\end{df}

\begin{df} 
We define local occupation number variables
\bel{lonv}
a_k(\bfeta) := \delta_{\eta(k),0}, \quad v_k(\bfeta) := \delta_{\eta(k),1}, \quad 
b_k(\bfeta) := \delta_{\eta(k),2}
\ee
and the global particle and vacancy numbers
\bel{partnum}
N(\bfeta) = \sum_{k=-L+1}^L a_k, \quad M(\bfeta) = \sum_{k=-L+1}^L b_k, \quad V(\bfeta) = \sum_{k=-L+1}^L v_k.
\ee
\end{df}
The argument of the local occupation number variables will be suppressed throughout this paper,
but not the argument of the global particle and vacancy numbers.
We note, for $\bfz=\bfeta$, the trivial but frequently used identities
\bel{partnum2}
N(\bfeta) \equiv N(\bfz) = |\bfx|, \quad M(\bfeta) \equiv M(\bfz) = |\bfy|,
\ee
\bel{occupos}
a_k = \sum_{i=1}^{N(\bfz)} \delta_{x_i,k} , \quad
 b_k  = \sum_{i=1}^{M(\bfz)} \delta_{y_i,k}.
\ee

\begin{df} 
For configurations $\bfz=\{\bfx,\bfy\}$ we define the number $N_{k}(\bfz)$ of 
$A$-particles to the left of site $k \in \Lambda$ and analogously
the number $M_{k}(\bfz)$ of $B$-particles to the left of site $k \in \Lambda$
\bel{NyMx}
N_{k}(\bfz) :=  \sum_{l=-L+1}^{k-1} a_l = \sum_{i=1}^{N(\bfz)} \sum_{l=-L+1}^{k-1} \delta_{x_i,l}, \quad 
M_{k}(\bfz) :=  \sum_{l=-L+1}^{k-1} b_l = \sum_{i=1}^{M(\bfz)} \sum_{l=-L+1}^{k-1} \delta_{y_i,l},
\ee
and
\bel{ABxyr}
A_k(\bfz) := 2 N_{k}(\bfz) - N(\bfz), \quad B_k(\bfz) := 2 M_{k}(\bfz) - M(\bfz).
\ee
\end{df}
Notice that the functions $N(\cdot), N_k(\cdot), a_k(\cdot), A_k(\cdot)$ depend only on the
$x$-coordinates (positions of the $A$-particles) of a configuration $\bfz$, while 
$M(\cdot), M_k(\cdot), b_k(\cdot), B_k(\cdot)$ depend only on the $y$-coordinates.

\subsection{Definition of the two-component ASEP}
\label{Sec:Defprocess}

Recalling the definitions \eref{etaexch} and \eref{occupos} the two-component ASEP $\bfeta_t$ 
described informally in the introduction is defined by the generator
\bel{2ASEPdef}
\mathcal{L} f(\bfeta)
= \sum_{k=-L+1}^{L-1} w^{kk+1}(\bfeta) [f(\bfeta^{kk+1}) - f(\bfeta)]
\ee
with the local hopping rates
\bel{localrates}
w^{kk+1}(\bfeta) = r \left( a_k v_{k+1} + v_k b_{k+1} + a_k b_{k+1} \right) +
\ell  \left(v_k a_{k+1} + b_k v_{k+1} + b_k a_{k+1} \right)
\ee
for a transition from a configuration 
$\bfeta$ to a configuration $\bfeta'$ with transition rate
\bel{transitionrate}
w(\bfeta\to\bfeta') = \sum_{k=-L+1}^{L-1} w^{kk+1}(\bfeta) \delta_{\bfeta',\bfeta^{kk+1}}.
\ee
It will turn out to be convenient to introduce the
asymmetry parameter $q$ and time-scale factor $w$
\bel{qw}
q = \sqrt{\frac{r}{\ell}}, \quad w = \sqrt{r\ell}.
\ee
The time scale will play no significant role below.

The general form of the evolution equation of a Markov chain with state space $\Omega$ and transition
rates $w(\bfeta\to\bfeta')$ for a transition from a configuration $\bfeta \in \Omega$ to a configuration
$\bfeta' \in \Omega$ is
\bel{generator}
\mathcal{L} f(\bfeta)
= {\sum}_{\eta' \in \Omega}' w(\eta\to\eta') [f(\eta') - f(\eta)]
\ee
where the prime at the summation indicates the absence of the term $\eta'=\eta$.
We define the transition matrix $H$ of the process by the matrix elements
\bel{transmatrix}
H_{\bfeta'\bfeta} = \left\{ \ba{ll} 
- w(\bfeta\to\bfeta') \quad & \bfeta \neq \bfeta' \\
{\sum}'_{\bfeta'} w(\bfeta\to\bfeta') & \bfeta = \bfeta' .
\ea \right.
\ee
The defining equation \eref{generator} then becomes 
\bel{generator2}
\mathcal{L} f(\bfeta) = - \sum_{\bfeta' \in \Omega} f(\bfeta') H_{\bfeta'\bfeta}
\ee
The r.h.s. of \eref{generator2} represents the multiplication of the matrix
$H$ with a vector whose components are $f(\bfeta')$ in the canonical basis.
In slight abuse of language we shall call also $H$ the generator of the process.
Below we shall construct $H$ for the two-component exclusion process 
in a judiciously chosen basis.

\subsection{The quantum algebra $U_q[\mathfrak{gl}(n)]$}

For the Lie algebra $\mathfrak{gl}(n)$ the quantum algebra
$U_q[\mathfrak{gl}(n)]$ is the associative algebra over $\C$
generated by $\mathbf{L}_i^{\pm 1}$, $i=1,\dots,n$ and $\mathbf{X}^\pm_i$, $i=1,\dots,n-1$
with the relations \cite{Jimb86,Burd92}
\bea
\label{Uqglndef1}
& & \comm{\mathbf{L}_i}{\mathbf{L}_j} = 0  \\
\label{Uqglncomm2}
& & \mathbf{L}_i \mathbf{X}^\pm_j = q^{\pm (\delta_{i,j+1} - \delta_{i,j})/2} \mathbf{X}^\pm_j \mathbf{L}_i\\
\label{Uqglncomm3}
& & \comm{\mathbf{X}^+_i}{\mathbf{X}^-_j} = 
\delta_{ij} \frac{(\mathbf{L}_{i+1}\mathbf{L}_i^{-1})^2 - (\mathbf{L}_{i+1}\mathbf{L}_i^{-1})^{-2}}{q-q^{-1}}
\eea
and, for $1 \leq i,j \leq n-1$,  the quadratic and cubic Serre relations 
\bea
\label{UqglnSerre1}
& & \comm{\mathbf{X}^\pm_i}{\mathbf{X}^\pm_j} = 0 \quad |i-j| \neq 1, \\
\label{UqglnSerre2}
& & (\mathbf{X}^\pm_i)^2 \mathbf{X}^\pm_j - [2]_q \mathbf{X}^\pm_i \mathbf{X}^\pm_j \mathbf{X}^\pm_i + 
\mathbf{X}^\pm_j (\mathbf{X}^\pm_i)^2 = 0 \quad |i-j| = 1.
\eea
Here the symmetric $q$-number is defined by
\bel{qnumber}
[x]_q := \frac{q^x - q^{-x}}{q-q^{-1}}
\ee
for $q,\,q^{-1} \neq 0$ and $x\in\C$. 
(Notice the replacement $q^2 \to q$ that we made in the definitions of \cite{Burd92}.)
The notion of symmetry of the generator under the action of the algebra means that 
there exist representation matrices $Y^\pm_i$ and $L_j$
of the algebra that all commute with the transition matrix $H$ of the process, i.e.,
\bel{Hsym}
\comm{H}{Y_i^\pm} = \comm{H}{L_i} = 0.
\ee
For the present case $n=3$ these representation matrices, given in
\eref{repY} and \eref{repL}, were constructed in \cite{Beli15a}.

\section{Results}

In order to state the first main result 
we first define the $q$-factorial
\bel{qfactorial}
[n]_q! := \left\{ \ba{ll} 1 & \quad n=0 \\
\prod_{k=1}^n [k]_q & \quad n \geq 1 \ea \right.
\ee
and the $q$-multinomial coefficients
\bel{qmultinomial}
C_K(N) = \frac{[K]_q!}{[N]_q![K-N]_q!}, \quad
C_{K}(N,M) = \frac{[K]_q!}{[N]_q![M]_q![K-N-M]_q!}.
\ee

\begin{theo}
\label{Theo:Invmeas}
The two-component exclusion process \eref{2ASEPdef} restricted to the subset $\mathbb{S}^{2L}_{N,M}$ 
of $N$ particles of type $A$ and $M$ particles of type $B$ has the unique invariant measure
\bel{theo1}
\pi^\ast_{N,M}(\bfeta) = \frac{\mathbf{1}_{\mathbb{S}^{2L}_{N,M}}(\bfeta)}{Z_{2L}(N,M)} \pi (\bfeta).
\ee
with the reversible measure
\be
\label{theo1b}
\pi(\bfeta) = q^{\sum_{k=-L+1}^L \left(2 k -1 \right)\left(a_k - b_k\right) 
+  \sum_{k=-L+1}^{L-1} \sum_{l=-L+1}^k \left( a_l b_{k+1}  -  b_l a_{k+1} \right)} 
\ee
and the normalization factor 
\bel{partfunNM}
Z_{2L}(N,M) = C_{2L}(N,M).
\ee
\end{theo}

We shall call these invariant measures, characterized in the following theorem, 
the canonical equilibrium distributions of the process.
Particle number conservation yields the following corollary.

\begin{cor}
\label{coro1}
The convex combinations
\bel{grandcan}
Q^\ast_{\nu,\mu} (\bfeta) =  \sum_{N=0}^{2L} \sum_{M=0}^{2L-N} 
 \frac{\rme^{\nu N + \mu M} Z_{2L}(N,M)}{Y_{2L}(\nu,\mu)}  \pi^\ast_{N,M}(\bfeta) =
\frac{\rme^{\nu N(\bfeta) + \mu M(\bfeta)}}{Y_{2L}(\nu,\mu)} \pi (\bfeta)
\ee
with the normalization factor
\bel{partfunnumu}
Y_{2L}(\nu,\mu) = \sum_{N=0}^{2L} \sum_{M=0}^{2L-N} \rme^{\nu N + \mu M} Z_{2L}(N,M)
\ee
are invariant measures for the two-component exclusion process \eref{2ASEPdef}.
\end{cor}

The second equality in \eref{grandcan} follows from the trivial identity
$\rme^{\nu N + \mu M} \mathbf{1}_{\mathbb{S}^{2L}_{N,M}}(\bfeta) = 
\rme^{\nu N(\bfeta) + \mu M(\bfeta)} \mathbf{1}_{\mathbb{S}^{2L}_{N,M}}(\bfeta)$.
We call these measures the grandcanonical equilibrium distributions. The normalization $Y_{2L}(\nu,\mu)$,
is a homogeneous bivariate Rogers-Szeg\H{o} polynomial \cite{Gasp04} and is called the
grandcanonical partition function. 

The limits $\mu\to -\infty$ or $\nu\to -\infty$ lead to the
pure grandcanonical measures
\bea
\label{pureA}
Q^{A \ast}_{\nu} (\bfeta) & = & \sum_{N=0}^{2L} 
\frac{\rme^{\nu N} C_{2L}(N)}{X_{2L}(\nu)}  \pi^\ast_{N,0}(\bfeta)\\
\label{pureB}
Q^{B \ast}_{\mu} (\bfeta) & = & \sum_{M=0}^{2L} 
\frac{\rme^{\mu M} C_{2L}(M)}{X_{2L}(\mu)}  \pi^\ast_{0,M}(\bfeta)
\eea
with the Rogers-Szeg\H{o} polynomial $X_{2L}(\alpha) = \sum_{K=0}^{2L} \rme^{\alpha K} C_{2L}(K)$.
From \eref{theo1} follows that $\pi^\ast_{N,0}(\bfeta) = 
\mathbf{1}_{\mathbb{S}^{2L}_{N,0}}(\bfeta) \tilde{\pi}_0(\bfeta)/Z_{2L}(N,0)$
with $\tilde{\pi}_0(\bfeta) = q^{\sum_{k=-L+1}^L \left(2 k -1 \right) a_k}$.
Since $Z_{2L}(N,0) = C_{2L}(N)$ one finds that $Q^{A \ast}_{\nu} (\bfeta)$
is a product measure in $\S^{2L}_{N,0}$ with marginals
$Q^{A \ast}_{\nu}(k) = (1+a_k(\rme^{\nu} q^{2k-1}-1))/(\rme^{\nu} q^{2k-1}+1)$, 
reminiscent of the blocking measure
of the single-species ASEP on $\Z$ \cite{Ligg85}. Likewise $Q^{B \ast}_{\mu} (\bfeta)$
is a product measure in $\S^{2L}_{0,M}$ with marginals
$Q^{B \ast}_{\mu}(k) = (1+b_k(\rme^{\mu} q^{-2k+1}-1))/(\rme^{\mu} q^{-2k+1}+1)$. 
The density profiles $\exval{a_k}_{\nu,0}$
and $\exval{b_k}_{0,\mu}$ in the 
pure grandcanonical measures follow by straightforward computation. One has shock
profiles
\bea
\exval{a_k}_{\nu,0} & = & \frac{\rme^{\nu} q^{2k-1}}{1+\rme^{\nu} q^{2k-1}} 
= \half \left[ 1 + \tanh{\left(\frac{k-\kappa_A}{\xi}\right)}\right] \\
\exval{b_k}_{0,\mu} & = & \frac{\rme^{\mu} q^{-2k+1}}{1+\rme^{\mu} q^{-2k+1}} 
= \half \left[ 1 - \tanh{\left(\frac{k-\kappa_B}{\xi}\right)}\right]
\eea
with the shock width $\xi = 1/\ln{q}$ and shock positions $\kappa_A = (1-\nu/\ln{q})/2$, 
$\kappa_B = (1+\nu/\ln{q})/2$.

In order to describe the self-duality of the process we define for configurations
$\bfeta \in \mathbb{S}^{2L}$ the functions
\bel{QAB}
Q_{x}^A(\bfeta) = q^{\sum_{k=-L+1}^{x-1} a_k - \sum_{k=x+1}^{L} a_k} a_x , \quad 
Q_{y}^B(\bfeta) = q^{-\sum_{k=-L+1}^{y-1} b_k + \sum_{k=y+1}^{L} b_k} b_y.
\ee
From these functions we construct the product
\bel{Q}
Q_{\bfz}(\bfeta) := \prod_{i=1}^{N(\bfz)} Q_{x_i}^A(\bfeta) \prod_{i=1}^{M(\bfz)} Q_{y_i}^B (\bfeta)
\ee
indexed by
$\bfz = \{\bfx,\bfy\}$, interpreted as a set of coordinates $x_i,y_i \in \Lambda$ and
unrelated to $\bfeta$.
With this definition we are in a position to state the second main result of this work.

\begin{theo}
\label{Theo:Duality}
Let $\bfz$ and $\bfeta$ be two configurations of the two-component exclusion 
process defined by 
\eref{2ASEPdef} with asymmetry parameter 
\eref{qw}. The process is self-dual with respect to the family of duality functions
\bel{theo2}
D(\bfz,\bfeta) = \pi^{-1}(\bfz) Q_{\bfz}(\bfeta)
\ee
where
$\pi^{-1}(\bfz)$ is the reversible measure \eref{theo1b}.
\end{theo}

We remark that the reversible measure \eref{theo1b} can be expressed as
\bel{theo1c}
\pi(\bfz) = q^{\sum_{i=1}^{N(\bfz)} [2x_i - 1 - M_{x_i}(\bfz)] - 
\sum_{i=1}^{M(\bfz)} [(2y_i -1 - N_{y_i}(\bfz)] }
\ee
by using \eref{occupos}.
Particle number conservation trivially induces independent duality relations for each
combination of particle number pairs $(N,M) = (N(\bfeta),M(\bfeta))$ and 
$(N',M')=(N(\bfz),M(\bfz))$ with duality functions 
\be
D_{N,M}^{N',M'}(\bfz,\bfeta):= 
D(\bfz_{N',M'},\bfeta_{N,M}) \mathbf{1}_{\S^{2L}_{N,M}}(\bfeta) \mathbf{1}_{\S^{2L}_{N',M'}}(\bfz). 
\ee
Therefore we refer to a ``family of duality functions''  rather than
just the ``duality function''. One has $D_{N,M}^{N',M'}(\bfz,\bfeta)=0$ if $N'>N$ or $M'>M$.
Using particle number conservation one can construct similar duality functions from
$\tilde{Q}^{A}_x := Q^{A}_x q^{N(\bfeta)} = q^{2\sum_{k=-L+1}^{x-1} a_k} a_x$ and 
$\tilde{Q}^{B}_y := Q^{B}_x q^{-M(\bfeta)} = q^{-2\sum_{k=-L+1}^{y-1} b_k} b_y$. For the
one-component ASEP $\tilde{Q}^{A}_x$ is the duality function of \cite{Schu97}.
A family of duality functions for the ASEP with second-class particles is given by
\eref{theo2} via the replacement $b_k \to v_k$. Duality and reversibility of the process
yield the following corollary, see \eref{coroproof}:

\begin{cor}
\label{Coro:Q}
For an initial distribution $P_0(\bfeta)$ with an arbitrary number 
$N$ of particles of type $A$ and $M$ particles of type $B$ we have for the time-dependent expectation 
\be
\label{coroQ}
\exval{Q_{\bfz}(t)}_{P_0}
= \sum_{z_{N,M}'} \exval{Q_{\bfz '}}_{P_0}
F(\bfz ;t| \bfz' ;0)
\ee
where $F(\bfz ;t| \bfz' ;0)$ is the transition probability of the two-component ASEP
with $N=N(\bfz)$ particles of type $A$ and $M=M(\bfz)$ particles of type $B$.
\end{cor}

Explicit exact expressions for $F(\bfz ;t| \bfz' ;0)$ have been obtained in \cite{Trac13} 
for the infinite system. 

Finally we present some simple properties characterizing the invariant measures. First we remark
that by the definition of the process -- in which the jumps of the $A$-particles ``do not see''
whether the neighbouring site is vacant or occupied by a $B$-particle -- 
expectations of the form $\exval{a_{k_1}\dots a_{k_n}}_{N,M}$ in the canonical 
equilibrium measure \eref{theo1} do not depend on $M$, i.e., 
$\exval{a_{k_1}\dots a_{k_n}}_{N,M} = \exval{a_{k_1}\dots a_{k_n}}_{N,0}$  Likewise, 
$\exval{b_{k_1}\dots b_{k_n}}_{N,M} = \exval{b_{k_1}\dots b_{k_n}}_{0,M}$ does not
depend on $N$.

The second characterization is a sum rule involving the the canonical invariant measures and the
duality function.

\begin{theo}
\label{Theo:Sumrule}
Let $\bfeta_{N,M}$ be a configuration in
$\mathbb{S}^{2L}_{N,M}$ with $N$ particles of type $A$ and
$M$ particles of type $B$ and let $\bfz$ be the coordinate representation of
a configuration in $\mathbb{S}^{2L}_{N',M'}$ with $N'$ particles of type $A$ and
$M'$ particles of type $B$. Then for all $\bfz \in \mathbb{S}^{2L}_{N',M'}$ and 
$\bfeta \in \mathbb{S}^{2L}_{N,M}$
one has the sum rule
\be
(\pi^\ast_{N',M'}(\bfz))^{-1} \sum_{\bfeta' \in \mathbb{S}^{2L}_{N,M}} 
\pi^\ast_{N,M}(\bfeta') Q_{\bfz}(\bfeta')
=  \sum_{\bfz' \in \mathbb{S}^{2L}_{N',M'}} Q_{\bfz'}(\bfeta) = \lambda_{N,M}^{N',M'} 
\ee
with a constant $\lambda_{N,M}^{N',M'}$ independent of $\bfeta$ and $\bfz$
and canonical stationary distribution given by \eref{theo1}.
\end{theo}

We remark that $\lambda_{N,M}^{N',M'}=0$ if $N'>N$ or $M'>M$.

\section{Some tools}

\subsection{More notation}

A generic time-dependent probability measure $\Prob{\bfeta_t = \bfeta}$ is denoted by
$P(\bfeta,t)$ or $P(\bfeta_t)$. For $t=0$ we use the notation $P_0(\bfeta):=P(\bfeta,0)$.
If $t$ is irrelevant we omit the argument $t$ and write $P(\bfeta)$. 
We also define the transition probability
\bel{transprob}
P(\bfeta',t|\bfeta,0)  := \Prob{\bfeta_t = \bfeta' | \bfeta_0 = \bfeta}
\ee
from a configuration $\eta$ to a configuration $\eta'$.

The expectation of a function $f(\bfeta)$ is denoted by
$\exval{f} := \sum_{\bfeta} f(\bfeta) P(\bfeta)$.
If we specify time and consider an initial distribution $P_0(\bfeta)$
we use for the expectation of a function $f(\bfeta_t)$ the notation
\bel{exval}
\exval{f(\bfeta_t)}_{P_0} := \sum_{\bfeta} P_0(\bfeta) \sum_{\bfeta'} f(\bfeta') P(\bfeta',t|\bfeta,0)
\ee 
or simply $\exval{f(t)}_{P_0}$.
For an initial distribution $P_0(\bfeta')= \delta_{\bfeta',\bfeta}$ concentrated on a configuration
$\bfeta$ we write $\exval{f(\eta_t)}_{\bfeta}$ or $\exval{f(t)}_{\bfeta}$.

\subsection{Matrix form of the generator}
\label{Sec:Matrix}

It turns out to be convenient to write the generator \eref{2ASEPdef} in the so-called quantum 
Hamiltonian form \cite{Schu01}
which is widely used in the physics literature on stochastic interacting particle systems and 
which was a 
given a formal probabilistic description in \cite{Lloy96}. However, this 
approach does not seem to be well-known in 
the probabilistic literature. For self-containedness and for introduction of our notation we summarize the
main ingredients.

\subsubsection{Choice of basis, inner product, and tensor product}

In order to write the matrix $H$ explicitly one has to choose an concrete basis,
i.e., to each configuration $\bfeta$ one has to assign a specific canonical basis vector.
Following \cite{Beli15a} we use ternary ordering, i.e., we assign to
each configuration $\bfeta$ the canonical basis vector 
\bel{ternary}
\iota(\bfeta) = 1+  \sum_{j=1}^{2L} \eta(j-L) 3^{j-1}.
\ee
of the complex
vector space $\C^{d}$ with dimension $d=3^L$. 
This basis vector has component 1 at position $\iota(\bfeta)$ and 0 else. We work with
a vector space over $\C$ rather than over $\R$ since in computations one encounters
eigenvectors of $H$ which may be complex.

We denote the basis vectors, which we consider to be column vectors, by $\ket{\bfeta}$.
We shall also use the notations $\ket{\bfz}$ and $\ket{\bfx,\bfy}$ instead of  $\ket{\bfeta}$. 
The basis vectors for configurations with a fixed number
$N$ of particles of type $A$ and $M$ particles of type $B$
are denoted by $\ket{\bfeta_{N,M}}$. 
We define also the dual basis $\bra{\bfeta} = \ket{\bfeta}^T$,
where the superscript $T$ on vectors or matrices denotes transposition.

The inner product of two arbitrary vectors
$\bra{w}$ with components $w_i$ and $\bra{v}$ with components $v_i$
is defined by
\bel{inprod}
\inprod{w}{v} = \sum_{i=1}^{d} w_i v_i
\ee
without complex conjugation.
In particular, we have the biorthogonality relation
\bel{inprodbasis}
\inprod{\bfeta}{\bfeta'} = \delta_{\bfeta \bfeta'}
\ee

Next we introduce the tensor product 
$\ket{v}\otimes \bra{w} \equiv \ket{v}\bra{w}$.
This tensor product is a $d \times d$-matrix with matrix elements 
$(\ket{v}\bra{w})_{i,j}=v_i w_j$.
Specifically we have the representation
\bel{unitmatrix}
\boldsymbol{1} = \sum_{\bfeta} \ket{\bfeta}\bra{\bfeta}.
\ee
of the $d$-dimensional unit matrix, expressing completeness of the basis.

\subsubsection{Measures and expectation values}

A probability measure $P(\bfeta)$ is represented by the probability vector
\bel{probvec}
\ket{P} = \sum_{\bfeta} P(\bfeta) \ket{\bfeta}.
\ee
From the inner product \eref{inprodbasis} and from \eref{generator2} we find
\bel{generator3}
\mathcal{L} f(\bfeta) = - \bra{f} H \ket{\bfeta}
\ee
where the vector $\bra{f} = \sum_{\bfeta} f(\bfeta) \bra{\bfeta}$ has components $f(\bfeta)$.
The semigroup property of the Markov chain is reflected in the time-evolution equation 
\be
\ket{P_t} = \rme^{-Ht} \ket{P_0} 
\ee
of a probability measure $P_0(\bfeta)$.

Normalization implies
\bel{normalization1}
\inprod{s}{P} = 1
\ee
where the {\it summation vector}
\bel{sumvec}
\bra{s} := \sum_{\bfeta} \bra{\bfeta}
\ee 
is the row vector where all components are equal to 1.
As a consequence one has
\bel{probcons}
\bra{s} H = 0
\ee
which means that the summation vector is a left eigenvector of $H$ with eigenvector 0.
This property follows from the fact that a diagonal element of $H_{\bfeta\bfeta}$ is by construction
the sum of all transition rates that appear with negative sign in the same column $\bfeta$ of $H$.
The vector corresponding to a stationary distribution is denoted $\ket{\pi^\ast}$. This is a right
eigenvector of $H$ with eigenvalue 0:
\bel{statvec}
H \ket{\pi^\ast} = 0.
\ee
and normalization $\inprod{s}{\pi^\ast}=1$.
An unnormalized right eigenvector with eigenvalue 0 is denoted $\ket{\pi}$.

The expectation $\exval{f}_P$ of a function $f(\bfeta)$ with respect to a probability distribution $P(\bfeta)$
becomes the inner product
\be
\exval{f}_P = \inprod{f}{P} = \bra{s} \hat{f} \ket{P}
\ee
where 
\bel{functionmatrix}
\hat{f} :=  \sum_{\bfeta} f(\bfeta)  \ket{\bfeta}  \bra{\bfeta}
\ee 
is a diagonal matrix with diagonal elements $f(\bfeta)$. Notice that
\bel{felement}
f(\bfeta) = \bra{\bfeta} \hat{f} \ket{\bfeta} = \bra{s}  \hat{f} \ket{\bfeta}.
\ee
For an initial distribution $P_0$ we can now use the definitions 
\eref{inprodbasis}, \eref{probvec}, \eref{sumvec} and the representation
\eref{unitmatrix} of the unit matrix to recover \eref{exval} in the matrix form
\bea
\exval{f(t)}_{P_0} & = & \sum_{\bfeta} P_0(\bfeta) 
\sum_{\bfeta'} f(\bfeta') P(\bfeta',t|\bfeta,0)  \nonumber \\
\label{exval2}
& = & \sum_{\bfeta'} \bra{s} \hat{f} \ket{\bfeta'} \bra{\bfeta'} \rme^{-Ht} \ket{\bfeta} \\ 
& = & \bra{s} \hat{f} \rme^{-Ht} \ket{P_0} \nonumber
\eea
Here 
\be
P(\bfeta',t|\bfeta,0) = \bra{\bfeta'} \rme^{-Ht} \ket{\bfeta}
\ee
is the transition probability \eref{transprob}. 

For a normalized stationary distribution we also define the diagonal matrix
\bel{invmeasmatrix}
\hat{\pi}^\ast :=  \sum_{\bfeta} \pi^\ast(\bfeta)  \ket{\bfeta}  \bra{\bfeta}.
\ee
For ergodic processes with finite state space one has $0< \pi^\ast(\bfeta) \leq 1$ for all $\bfeta$.
Then all powers $(\hat{\pi}^\ast)^\alpha$ exist.
In terms of this diagonal matrix we can write the generator of the reversed dynamics as 
\bel{revtrans2}
H^{rev}  = \hat{\pi}^\ast H^T (\hat{\pi}^\ast)^{-1}.
\ee
Reversibility means $H^{rev} = H$. An unnormalized stationary distribution $\pi$ for which 
\bel{reversibility}
H \hat{\pi} = \hat{\pi} H^T
\ee 
holds with
\bel{reversmatrix}
\hat{\pi} = \sum_{\bfeta} \pi(\bfeta)  \ket{\bfeta}  \bra{\bfeta}.
\ee
is called a reversible measure.

\subsubsection{Explicit form of the generator}

In order to write the generator $H$ explicitly we define
 the following matrices:
\bea
\label{creation}
& & a^+ := \left( \ba{ccc} 0 & 1 & 0 \\ 0 & 0 & 0 \\ 0 & 0 & 0 \ea \right), \quad
b^+ := \left( \ba{ccc} 0 & 0 & 0 \\ 0 & 0 & 0 \\ 0 & 1 & 0 \ea \right), \quad
c^+ := \left( \ba{ccc} 0 & 0 & 1 \\ 0 & 0 & 0 \\ 0 & 0 & 0 \ea \right), \\
\label{annihilation}
& & a^- := \left( \ba{ccc} 0 & 0 & 0 \\ 1 & 0 & 0 \\ 0 & 0 & 0 \ea \right), \quad
b^- := \left( \ba{ccc} 0 & 0 & 0 \\ 0 & 0 & 1 \\ 0 & 0 & 0 \ea \right),  \quad
c^- := \left( \ba{ccc} 0 & 0 & 0 \\ 0 & 0 & 0 \\ 1 & 0 & 0 \ea \right),
\eea
the diagonal projectors
\bea
\label{projection}
& & \hat{a} :=  \left( \ba{ccc} 1 & 0 & 0 \\ 0 & 0 & 0 \\ 0 & 0 & 0 \ea \right), \quad
\hat{v} : = \left( \ba{ccc} 0 & 0 & 0 \\ 0 & 1 & 0 \\ 0 & 0 & 0 \ea \right), \quad
\hat{b} :=  \left( \ba{ccc} 0 & 0 & 0 \\ 0 & 0 & 0 \\ 0 & 0 & 1 \ea \right). 
\eea
and the three-dimensional unit matrix
\bel{3dunit}
\mathds{1} = \hat{a} + \hat{v} + \hat{b}.
\ee

For matrices $M$ the expression $M^{\otimes j}$ will denote the $j$-fold tensor
product of $M$ with itself if $j>1$. For $j=1$ we define $M^{\otimes 1} := M$ and for $j=0$ we define 
$M^{\otimes 0} = 1$
with the $c$-number $1$.
For arbitrary $3\times 3$-matrices $u$
we define tensor operators
\be
u_k := \mathds{1}^{\otimes (k+L-1)} \otimes u \otimes \mathds{1}^{\otimes (L-k)}
\ee
which allow us to write the generator $H$ for the two-component ASEP on the lattice $\{-L+1,\dots,L\}$ 
as \cite{Beli15a}
\bel{2ASEPgen}
H = \sum_{k=-L+1}^{L-1} h_{k,k+1}
\ee
with the hopping matrices
\bea
h_{k,k+1} & : = & r \left( \hat{a}_k \hat{v}_{k+1} - a^-_k a^+_{k+1}  + 
\hat{v}_k \hat{b}_{k+1} - b^+_k b^-_{k+1} 
+ \hat{a}_k \hat{b}_{k+1} - c^-_k c^+_{k+1}  \right) \nonumber \\
\label{hoppingmatrixk}
& & +
\ell  \left(\hat{v}_k \hat{a}_{k+1} - a^+_k a^-_{k+1}  + \hat{b}_k \hat{v}_{k+1} - b^-_k b^+_{k+1} 
+  \hat{b}_k \hat{a}_{k+1} - c^+_k c^-_{k+1}  \right).
\eea

With \eref{qw} we split $H = H_d + H_o$ into its offdiagonal part 
\bel{Hoff}
H_o = - w \sum_{k=-L+1}^{L-1} \left[ q \left(a^-_k a^+_{k+1}  + b^+_k b^-_{k+1} + c^-_k c^+_{k+1}  \right)
+ q^{-1}  \left(a^+_k a^-_{k+1}  + b^-_k b^+_{k+1} +  c^+_k c^-_{k+1}  \right) \right]
\ee
and its diagonal part
\bel{Hdiag}
H_d =  w \sum_{k=-L+1}^{L-1} \left[ q \left( \hat{a}_k \hat{v}_{k+1}  + \hat{v}_k \hat{b}_{k+1} 
+ \hat{a}_k \hat{b}_{k+1}  \right)
+ q^{-1}  \left( \hat{v}_k \hat{a}_{k+1}  + \hat{b}_k \hat{v}_{k+1} 
+  \hat{b}_k \hat{a}_{k+1}  \right) \right].
\ee

For more details of the construction of $H$ in the tensor basis we refer the reader to \cite{Beli15a}.

\subsection{Duality}
\label{Sec:Duality}

We recall the concept of duality in matrix form \cite{Sudb95,Giar09}, see also \cite{Jans14} for
a detailed discussion.
In this subsection $X$ and $\Omega$ represent arbitrary finite-dimensional state spaces.
Consider two processes $x_t$ and $\omega_t$ and a 
function $D: X \times \Omega \mapsto \C$. Notice that the function $D(x,\omega)$ can be
understood as a family of functions $f_x:\Omega \mapsto \C$ indexed by $x$ and defined by 
$f_x(\omega) := D(x,\omega)$, or,
alternatively as a family of functions $g_\omega:X \mapsto \C$ indexed by $\omega$
and defined by $g_\omega(x) := D(x,\omega)$.

The two processes are said to be dual to each other if
\bel{Def:duality}
\exval{D(x,\omega_t)}_\omega = \exval{D(x_t,\omega)}_x
\ee
We remark that with the definitions introduced above we have
\bea
& & \exval{D(x,\omega_t)}_\omega = \sum_{\omega'} D(x,\omega') P(\omega',t|\omega,0) = \exval{f_x(t)}_\omega\\
& & \exval{D(x_t,\omega)}_x = \sum_{x'} D(x',\omega) P(x',t|x,0) = \exval{g_\omega(t)}_x
\eea
so that duality can be stated as
\bel{Def:duality2}
\exval{f_x(t)}_\omega = \exval{g_\omega(t)}_x
\ee
with $\exval{f_x(0)}_\omega = \exval{g_\omega(0)}_x = D(x,\omega)$.

In order to make contact with the quantum Hamiltonian formalism we
define $\ket{x}$ as a canonical basis vector of $\C^{|X|}$ and $\ket{\omega}$ as a canonical 
basis vector of $\C^{|\Omega|}$. Let $\bra{s}$ and $\bra{\tilde{s}}$ be the corresponding
summation vectors.
Define the matrix
\bel{dualitymatrix}
D = \sum_x \sum_\omega D(x,\eta) \ket{x}\bra{\omega}
\ee
with matrix elements $\bra{x} D \ket{\omega} = D(x,\omega)$.
The processes $\omega_t$  and $x_t$ with generators $H$ and $G$ 
are dual to each other w.r.t. the duality function
$D(x,\omega)$ if
\bel{Def:duality3}
DH = G^TD.
\ee

It is easy to prove the equivalence of this definition with the original definition
\eref{Def:duality}. Since the kind of arguments underlying this equivalence are
important for the present matrix formulation of duality we present them here in detail:
\bea
\exval{D(x,\omega_t)}_\omega = &  \sum_{\omega'} D(x,\omega') P(\omega',t|\omega,0) & \\
\label{step1}
= & \sum_{\omega'} \bra{x} D \ket{\omega'} \bra{\omega'} \rme^{-Ht} \ket{\omega} & \\
\label{step2}
= & \bra{x} D \rme^{-Ht} \ket{\omega} & \\
\label{step3}
= & \bra{x}  \rme^{-G^Tt} D \ket{\omega} & \\
\label{step4}
= & \sum_{x'} \bra{x} \rme^{-G^Tt} \ket{x'} \bra{x'} D\ket{\omega} & \\
\label{step5}
= & \sum_{x'} D(x',\omega) P(x',t|x,0) & = \exval{D(x_t,\omega)}_x
\eea
In going from \eref{step1} to \eref{step2} and from \eref{step4} to \eref{step5} we use 
the representation of the unit matrix constructed in analogy to \eref{unitmatrix}.
In the step from from \eref{step2} to \eref{step3} we apply the definition \eref{Def:duality3}.
Since we have a chain of equalities, it can be read in both directions. Thus the
equivalence is established. 

In order to express the alternative definition \eref{Def:duality2} in matrix form
we introduce the diagonal matrices
\bel{dualdiag}
\hat{f}_x = \sum_\omega D(x,\omega) \ket{\omega}\bra{\omega}, \quad
\hat{g}_\omega = \sum_x D(x,\omega) \ket{x}\bra{x}.
\ee
The duality relation \eref{Def:duality2} reads
\bel{Def:duality4}
\bra{s} \hat{f}_x \rme^{-Ht} \ket{\omega} = \bra{\tilde{s}} \hat{g}_\omega \rme^{-Gt} \ket{x}.
\ee

To prove of equivalence of \eref{Def:duality4} with \eref{Def:duality3} we
note that by construction
\bel{operatorduality}
\bra{s} \hat{f}_x = \bra{x} D, \quad \hat{g}_\omega = \ket{\tilde{s}} = D \ket{\omega}.
\ee
Then it follows that
\bea 
\bra{s} \hat{f}_x \rme^{-Ht} \ket{\omega} & \bra{x} D \rme^{-Ht} \ket{\omega} & \\
= & \bra{x}  \rme^{-G^Tt} D \ket{\omega} & \\
= & \bra{x} \rme^{-G^Tt} \hat{g}_\omega  \ket{\tilde{s}} & = 
\bra{\tilde{s}} \hat{g}_\omega \rme^{-Gt} \ket{x}
\eea
which establishes the equivalence.

We end this discussion with a reformulation of Theorem 2.6 of \cite{Giar09}.

\begin{theo}
\label{Theo:Dualrev}
Let $H$ be the matrix representation of the generator of an ergodic Markov process $\eta_t$ 
with countable state space and $H^{rev}$
be the matrix form of the generator of the reversed process $\xi_t$. 
Assume that there exists an intertwiner $S$
such that 
\bel{hypoS}
S H = H^{rev} S.
\ee
Then $H$ is self-dual with duality function $D(\xi,\eta) = D_{\xi,\eta}$ given by the
matrix elements of the duality matrix
\bel{Theodualrevmatrix}
D = \hat{\pi}^{-1} S.
\ee
with the diagonal stationary distribution matrix \eref{reversmatrix}.
\end{theo}

The proof that $SH=H^{rev}S$ implies self-duality with duality matrix $D=(\hat{\pi}^\ast)^{-1}S$
is elementary and follows from the chain of equalities
\be 
DH = \hat{\pi}^{-1} S H = \hat{\pi}^{-1} H^{rev} S = 
\hat{\pi}^{-1} H^{rev} \hat{\pi} D = H^T D.
\ee
The first and the third equality are the definition 
\eref{Theodualrevmatrix}, the second equality is the hypothesis
\eref{hypoS} of the theorem, and the fourth equality is the reversibility relation \eref{reversibility}.

\begin{rem}
It follows that if $H$ is reversible then the hypothesis \eref{hypoS} reads $SH=HS$, i.e. 
$S$ is a symmetry of $H$.
Unlike \cite{Giar09} we do not require $S$ to be invertible. 
\end{rem}

\subsection{Representation matrices for $U_q[\mathfrak{gl}(3)]$}

\subsubsection{Relation between $U_q[\mathfrak{gl}(n)]$ and $U_q[\mathfrak{sl}(n)]$}

It is convenient to introduce generators $\mathbf{H}_i$ and $\tilde{\mathbf{H}}_i$ through
\bel{complement}
q^{-\tilde{\mathbf{H}}_i/2} = \mathbf{L}_i, \quad \mathbf{H}_i = 
\tilde{\mathbf{H}}_i - \tilde{\mathbf{H}}_{i+1}.
\ee
Then the quantum algebra $U_q[\mathfrak{sl}(n)]$ is the subalgebra generated by 
$q^{\pm \mathbf{H}_i/2}$, and $\mathbf{X}^\pm_i$, $i=1,\dots,n-1$ with relations
\eref{UqglnSerre1}, \eref{UqglnSerre2}
and
\bea
& & q^{\mathbf{H}_i/2} q^{-\mathbf{H}_i/2} = q^{-\mathbf{H}_i/2} q^{\mathbf{H}_i/2} = I\\
\label{Uqslncomm1b}
& & q^{\mathbf{H}_i/2} q^{\mathbf{H}_j/2} = q^{\mathbf{H}_j/2} q^{\mathbf{H}_i/2} \\
\label{Uqslncomm2b}
& & q^{\mathbf{H}_i} \mathbf{X}^\pm_j q^{-\mathbf{H}_i} = q^{\pm A_{ij}} \mathbf{X}^\pm_j\\
\label{Uqslncomm3b}
& & \comm{\mathbf{X}^+_i}{\mathbf{X}^-_j} = \delta_{ij} [\mathbf{H}_i]_q.
\eea
with the unit $I$ and the Cartan matrix $A$ of simple Lie algebras of type $A_n$
\bel{basicdefs}
 A_{ij} :=  \left\{ \ba{rl} 2 & i=j \\ -1 & j = i\pm 1\\ 0 & \mbox{else.} \ea \right.
\ee
The fact that $U_q[\mathfrak{sl}(n)]$ is a subalgebra of $U_q[\mathfrak{gl}(n)]$ can be seen by noticing that 
$\sum_{i=1}^n \tilde{\mathbf{H}}_i$ belongs to the center of $U_q[\mathfrak{gl}_n]$ \cite{Jimb86}.

\subsubsection{Tensor representation for $n=3$}

In order to distinguish the three-dimensional matrices corresponding to  the
fundamental representation 
from the abstract generators we use lower case letters. 
In terms of \eref{creation}, \eref{annihilation}, \eref{projection}
the fundamental representation of $U_q[\mathfrak{gl}(3)]$ is given by:
\bea
\label{fundrepX}
& & x_1^\pm = a^\pm, \quad x_2^\pm = b^\mp \\
\label{fundrepH3}
& & \tilde{h}_1 = \hat{a}, \quad \tilde{h}_2 = \hat{v}, \quad \tilde{h}_3 = \hat{b}, 
\eea
corresponding to
\be 
\label{fundrepH}
h_1 = \hat{a} - \hat{v}, \quad
h_2 = \hat{v} - \hat{b}.
\ee
for the representation of the generators $\mathbf{H}_i$ of $U_q[\mathfrak{sl}(3)]$.
It is convenient to work both with $h_{i}$ and the projectors $\tilde{h}_i$
expressed in term of the projectors \eref{projection}.

In terms of the fundamental representation a tensor representation of 
$U_q[\mathfrak{sl}(3)]$, denoted by boldface capital letters,
is given by \cite{Beli15a}
\bel{repY}
Y^\pm_i = \sum_{k=-L+1}^L  Y^\pm_i(k)
\ee
with
\bea
\label{X1ptrafo}
Y^+_1(k)  & = & q^{ \sum_{l=-L+1}^{k-1} \hat{v}_l
- \sum_{l=k+1}^{L} \hat{v}_l} a_k^+, \\
\label{X1mtrafo}
Y^-_1(k)  & = & q^{- \sum_{l=-L+1}^{k-1} \hat{a}_l
+ \sum_{l=k+1}^{L} \hat{a}_l} a_k^-, \\
\label{X2ptrafo}
Y^+_2(k)  & = & q^{ \sum_{l=-L+1}^{k-1} \hat{b}_l 
- \sum_{l=k+1}^{L} \hat{b}_l} b_k^- \\
\label{X2mtrafo}
Y^-_2(k)  & = & q^{- \sum_{l=-L+1}^{k-1} \hat{v}_l 
+ \sum_{l=k+1}^{L} \hat{v}_l} b_k^+
\eea
and
\bel{repH}
H_i = \sum_{k=-L+1}^L  H_i(k)
\ee
with
\be
H_i(k) = \mathds{1}^{\otimes k+L-1} \otimes h_i \otimes \mathds{1}^{ \otimes L-k}.
\ee
Notice that $H_1(k) = \hat{a}_k - \hat{v}_k$ and $H_2(k) = \hat{v}_k - \hat{b}_k$.

For the full quantum algebra $U_q[\mathfrak{gl}(3)]$ we have the diagonal representation matrices
\bel{repH3}
\tilde{H}_1  = \sum_{k=-L+1}^L \hat{a}_k =: \hat{N}, 
\quad \tilde{H}_2 = \sum_{k=-L+1}^L \hat{v}_k, 
\quad \tilde{H}_3 = \sum_{k=-L+1}^L \hat{b}_k =: \hat{M}.
\ee
Here $\hat{N}$ and $\hat{M}$ are the particle number operators satisfying
\be 
\hat{N} \ket{\bfeta_{N,M}} = N \ket{\bfeta_{N,M}}, \quad \hat{M} \ket{\bfeta_{N,M}} = M \ket{\bfeta_{N,M}}.
\ee
From these matrices one obtains the representation matrices
\bel{repL}
L_i=q^{-\tilde{H}_i/2}.
\ee
The unit $I$ is represented by the $3^L$-dimensional unit matrix $\mathbf{1} := \mathds{1}^{\otimes 2L}$.

The crucial property of the representation \eref{repY} and \eref{repL}
that was proved in \cite{Beli15a} and which is used heavily below
are the commutation relations \eref{Hsym}
which express the symmetry of the generator $H$ \eref{2ASEPgen} under the action of the quantum algebra 
$U_q[\mathfrak{gl}(3)]$.

\section{Proofs}

\subsection{Proof of Theorem \eref{Theo:Invmeas}}

(i) We first note that uniqueness of $\pi^\ast_{N,M}$ follows from ergodicity of the 
process defined on the subset $\S^{2L}_{N,M}$ which is ensured by the fact that
the process is a random sequence of permutations $\sigma^{k,k+1}(\bfeta)$.

(ii) In \cite{Beli15a} we proved, using the quantum algebra symmetry \eref{Hsym}, 
that the two-component exclusion process defined by 
\eref{2ASEPdef} has the unnormalized reversible measure $\pi$ \eref{theo1b}.
Below we give a direct proof without reference to the quantum algebra symmetry.
According to the discussion of Section \eref{Sec:Matrix}
we prove the transformation property \eref{reversibility} with the 
generator \eref{2ASEPgen}. Since $\hat{\pi}$ is diagonal one has
has $\hat{\pi}^{-1} H_d \hat{\pi} = H_d$ for the diagonal part
\eref{Hdiag} of $H$. It remains to show that $\hat{\pi}^{-1} H_o \hat{\pi} = H_o^T$
for the offdiagonal part \eref{Hoff}.
To this end we first prove the basic transformation lemma

\begin{lmm}
\label{Lemma:abtrafolmx}
For any finite $p\neq 0$ we have 
\bea
\label{abtrafo1a}
& & p^{\hat{a}_l} a_x^\pm p^{-\hat{a}_l} = p^{\pm \delta_{l,x}} a_x^\pm, \quad
p^{\hat{b}_l} a_x^\pm p^{-\hat{b}_l} =  a_x^\pm, \\
\label{abtrafo1b}
& & p^{\hat{b}_l} b_x^\pm p^{-\hat{b}_l} = p^{\pm \delta_{l,x} } b_x^\pm, \quad
p^{\hat{a}_l} b_x^\pm p^{-\hat{a}_l} =  b_x^\pm \\
\label{abtrafo2a}
& & p^{\hat{a}_l \hat{b}_m} a_x^\pm p^{-\hat{a}_l \hat{b}_m} = p^{\pm \delta_{l,x} \hat{b}_m} a_x^\pm,\\
\label{abtrafo2b}
& & p^{\hat{a}_l \hat{b}_m} b_x^\pm p^{-\hat{a}_l \hat{b}_m} = p^{\pm \delta_{m,x} \hat{a}_l} b_x^\pm.
\eea
\end{lmm}

{\it Proof:}
A projector $\hat{u}$ has the property $\hat{u}=\hat{u}^2$.
Thus its exponential can be written $p^{\hat{u}} = 1 + (p-1)\hat{u}$.
Since $\hat{a}_l \hat{b}_{m}$ is a projector
one has 
\bel{expproj}
p^{\hat{a}_l \hat{b}_{k+1}} = 1 + (p-1)\hat{a}_l \hat{b}_{k+1}.
\ee
The tensor construction implies that $u_k u_l = u_l u_k$ for $k\neq l$ and any $u$. 
For $k=l$ we observe that  one obtains by direct computation the relations
\bea
\label{projhatar}
a^+ \hat{a} = b^+ \hat{a} = b^- \hat{a} = c^+ \hat{a} = 0, &\quad a^- \hat{a} = a^-, 
&\quad c^- \hat{a} = c^-\\
a^- \hat{\upsilon} = b^- \hat{\upsilon} = c^+ \hat{\upsilon} = c^- \hat{\upsilon} = 0, 
&\quad a^+ \hat{\upsilon} = a^+, &\quad b^+ \hat{\upsilon} = b^+\\
\label{projhatbr}
a^+ \hat{b} = a^- \hat{b} = b^+ \hat{b} = c^- \hat{b} = 0, &\quad b^- \hat{b} = b^-, 
&\quad c^+ \hat{b} = c^+
\eea
and
\bea
\label{projhatal}
\hat{a} a^- = \hat{a} b^+ = \hat{a} b^- = \hat{a} c^- = 0, &\quad \hat{a} a^+ = a^+, 
&\quad \hat{a} c^+ =  c^+\\
\hat{\upsilon} a^+ = \hat{a} b^+ = \hat{\upsilon} c^+ = \hat{\upsilon} c^- = 0, 
&\quad \hat{\upsilon} a^- = a^-, &\quad \hat{\upsilon} b^- =  b^-\\
\label{projhatbl}
\hat{b} a^+ = \hat{b} a^- = \hat{b} b^- = \hat{b} c^+ = 0, &\quad \hat{b} b^+ = b^+, 
&\quad \hat{b} c^- =  c^-.
\eea
By multilinearity of the tensor product these relations remain valid on each subspace $k$.
Relations \eref{abtrafo1a} - \eref{abtrafo2b} then follow from \eref{expproj}. \qed

Now we decompose
\bel{pideco}
\hat{\pi} = \hat{A} \hat{B} \hat{U}
\ee
with
\be
\hat{A} =  q^{ \sum_{k=-L+1}^L \left(2 k -1 \right) a_k}, \,\,
\hat{B} =  q^{-\sum_{k=-L+1}^L \left(2 k -1 \right) b_k}, \,\,
\hat{U} = \prod_{k=-L+1}^{L-1} \prod_{l=1}^k q^{\hat{a}_l \hat{b}_{k+1} - \hat{b}_l \hat{a}_{k+1}}.
\ee

Together with $c^\pm = a^\pm b^\mp$ one has from 
\eref{abtrafo1a} and \eref{abtrafo1b} of Lemma \eref{Lemma:abtrafolmx} 
for $-L+1 \leq k \leq L$ 
\bea
\hat{A} a_k^\pm \hat{A}^{-1} = q^{\pm(2k-1)} a_k^\pm, & \hat{A} b_k^\pm \hat{A}^{-1} = b_k^\pm, &
\hat{A} c_k^\pm \hat{A}^{-1} = q^{\pm(2k-1)} c_k^\pm, \\
\hat{B} b_k^\pm \hat{B}^{-1} = q^{\mp(2k-1)} b_k^\pm, & \hat{B}  a_k^\pm \hat{B}^{-1} = a_k^\pm, &
\hat{B} c_k^\pm \hat{B}^{-1} = q^{\pm(2k-1)} c_k^\pm.
\eea
For the transformation $U_p$ one obtains from \eref{abtrafo2a} and \eref{abtrafo2b} of 
Lemma \eref{Lemma:abtrafolmx} for $-L+1 \leq k \leq L$ 
\bea
\hat{U} a_k^\pm \hat{U}^{-1} & = & 
q^{\mp \sum_{l=-L+1}^{k-1} \hat{b}_l  \pm \sum_{l=k+1}^{L} \hat{b}_l} a_k^\pm \\
\hat{U} b_k^\pm \hat{U}^{-1} & = & 
q^{\mp \sum_{l=-L+1}^{k-1} \hat{a}_l  \pm \sum_{l=k+1}^{L} \hat{a}_l} b_k^\pm \\
\hat{U} c_k^\pm \hat{U}^{-1} & = & 
q^{\mp \sum_{l=-L+1}^{k-1} (\hat{b}_l-\hat{a}_l)\pm \sum_{l=k+1}^{L} (\hat{b}_l-\hat{a}_l)} c_k^\pm .
\eea

Putting these results together and using the projector property \eref{expproj}
together with \eref{projhatar} - \eref{projhatbl} yield for $-L+1 \leq k \leq L-1$ 
\be
\hat{A} a_k^\pm a_{k+1}^\mp \hat{A}^{-1} = q^{\mp 2} a_k^\pm a_{k+1}^\mp, \,
\hat{A} b_k^\pm b_{k+1}^\mp \hat{A}^{-1} = b_k^\pm b_{k+1}^\mp, \,
\hat{A} c_k^\pm c_{k+1}^\mp \hat{A}^{-1} = q^{\mp 2} c_k^\pm c_{k+1}^\mp,
\ee
\be
\hat{B} b_k^\pm b_{k+1}^\mp \hat{B}^{-1} = q^{\pm 2} b_k^\pm b_{k+1}^\mp, \,
\hat{B} a_k^\pm a_{k+1}^\mp \hat{B}^{-1} = a_k^\pm a_{k+1}^\mp, \,
\hat{B} c_k^\pm c_{k+1}^\mp \hat{B}^{-1} = q^{\mp 2} c_k^\pm c_{k+1}^\mp
\ee
and
\be
\hat{U} a_k^\pm a_{k+1}^\mp \hat{U}^{-1} = a_k^\pm a_{k+1}^\mp,
\hat{U} b_k^\pm b_{k+1}^\mp \hat{U}^{-1} = b_k^\pm b_{k+1}^\mp ,
\hat{U} c_k^\pm c_{k+1}^\mp \hat{U}^{-1} = q^{\pm 2} c_k^\pm c_{k+1}^\mp. 
\ee
Since $(a^\pm)^T=a^\mp$ (and similarly for $b^\pm$ and $c^\pm$) applying the decomposition 
\eref{pideco} to the individual terms in \eref{2ASEPgen} yields
$\hat{\pi}^{-1} H_o \hat{\pi} = H_o^T$ and therefore reversibility of $\pi$.

(iii): We complete the proof of Theorem \eref{Theo:Invmeas} by proving the normalization
factor. For a configuration $\bfz=\{\bfx,\bfy\}$ 
define $\tilde{y}_i = = y_i - N_{y_i}(\bfz)$. 
We have by definition of the partition function
\be
Z_{2L}(N,M) = \sum_{\bfz_{N,M}} \pi(\bfz_{N,M}) =
\sum_{\bfz_{N,M}} q^{ \sum_{i=1}^N (2x_i-1)- \sum_{i=1}^M (2\tilde{y}_i +N-1)}.
\ee
Consider the points $\vec{r}$ in
the Weyl alcove $W_K^{2L} = \{\vec{r} : -L < x_1 < ... < x_K \leq L\}$. We also define the
punctuated Weyl alcove $W_K^{2L}(\vec{r}) = W_K^{2L}\setminus \vec{r}$ for $\vec{r} \in W_K^{2L}$.
This allows us to write 
$\sum_{\bfz_{N,M}}= \sum_{\vec{x}\in W_N^{2L}} \sum_{\vec{y}\in W_{M}^{2L}(\vec{x})}$.

Next observe that
by construction $\sum_{\vec{y}\in W_{M}^{2L}(\vec{x})} f(\tilde{y}_i) = \sum_{\vec{y}\in W_{M}^{2L-N}} f(y_i)$.
Therefore 
\be
Z_{2L}(N,M) = \sum_{\vec{x}\in W_N^{2L}} \sum_{\vec{y}\in W_{M}^{2L-N}} q^{ \sum_{i=1}^N (2x_i-1)- \sum_{i=1}^M (2y_i+N-1)} 
\ee
which implies that 
$Z_{2L}(N,M) = Z_{2L}(N,0) Z_{2L-N}(0,M)$. A classical result from the theory of integer 
partitions \cite{Andr76}
yields for the single-species partition functions $Z_{2L}(N,0) = C_{2L}(N)$, 
$Z_{2L}(0,M) C_{2L}(M)$ 
with the $q$-binomial coefficient $C_K(N)$.
Observing that $C_{2L}(N) C_{2L-N}(M) = C_{2L}(N,M)$ concludes the proof. \qed

\subsection{Proof of Theorem \eref{Theo:Duality}}

\subsubsection{Reformulation of the problem}

Step 1: We first apply the general considerations of Sec. \eref{Sec:Duality}
to the present case of the two-component
ASEP. It is convenient to use the occupation
variable presentation $\bfeta_t$ for one process and the coordinate representation $\bfz_t$
for the dual.
The duality function, given by $\bra{\bfz} D \ket{\bfeta}$ in terms of the duality matrix $D$, 
is therefore denoted by $D(\bfz,\bfeta)$.
If self-duality is valid
for some duality function $D(\bfz,\bfeta)$ then according to \eref{dualdiag}
we can define a diagonal matrix $\hat{D}_\bfz$ such that 
\bel{dualdiag1}
D(\bfz,\bfeta) = \bra{s} \hat{D}_\bfz \ket{\bfeta}
\ee
with the summation vector $\bra{s}$. Then self-duality yields
\bel{coroproof}
\bra{s}  \hat{D}_\bfz \rme^{-Ht} \ket{\bfeta}= \sum_{\bfz'} \bra{\bfz'}  
\rme^{-Ht} \ket{\bfz} \bra{s}  \hat{D}_\bfz  \ket{\bfeta}
\ee
and, as a consequence from reversibility, Corollary \eref{coroQ}.

Step 2: In \cite{Beli15a} we have established the symmetry of the generator under the action of
$U_q[\mathfrak{gl}(3)]$. Moreover, we have reversibility
$H=H^{rev}$ of the two-component ASEP  with the reversible measure \eref{theo1b}. 
Then for any matrix $S$ satisfying $[S,H]=0$ Theorem \eref{Theo:Dualrev} yields a duality
function
\be
D(\bfz,\bfeta) =
\bra{\bfz} (\hat{\pi}^\ast(\bfz) )^{-1} S \ket{\bfeta} = \pi^{-1}(\bfz) \bra{\bfz}  S \ket{\bfeta}.
\ee
which means that we can construct duality functions from the
symmetry operators of the model, i.e., from the tensor representation \eref{repY}, \eref{repH}.

Step 3: On the other hand, from \eref{dualdiag1}, one has $D(\bfz,\bfeta) = 
\bra{s} \hat{D}_\bfz \ket{\bfeta}$ for all $\bfeta \in \S^{2L}$. 
Therefore one can express the duality function 
\be
D(\bfz,\bfeta) = \pi^{-1}(\bfz) Q_\bfz(\bfeta) =
\pi^{-1}(\bfz)
\bra{s} \hat{Q}_\bfz \ket{\bfeta}. 
\ee
of Theorem \eref{Theo:Duality}
in terms of a diagonal matrix $\hat{Q}_\bfz$ satisfying
\bel{dualityalter}
\bra{\bfz} S = \bra{s} \hat{Q}_\bfz 
\ee
and $\bra{s} \hat{Q}_\bfz \ket{\bfeta} = Q_\bfz(\bfeta)$ given in \eref{Q}.
Therefore the task at hand is to find a symmetry operator $S$ that satisfies 
\eref{dualityalter} with the diagonal matrix $\hat{Q}_\bfz$ with matrix elements
given by \eref{Q}.

Step 4: In order to choose $S$ we observe that $D(\emptyset,\eta) = 1$,
corresponding to $\hat{D}_\emptyset = \hat{Q}_\emptyset = \mathds{1}$. 
The non-trivial information one gains is that
$ \bra{\emptyset} S =\bra{s}$ which means 
that the symmetry operator $S$ generates
the summation vector  from the vacuum vector $\bra{\emptyset}$. From the explicit representation
obtained in \cite{Beli15a} we find as a candidate
\bel{S}
S = \sum_{n=0}^{2L} \sum_{m=0}^{2L-n} \frac{(Y^-_1)^n}{[n]_q!} \frac{(Y^+_2)^m}{[m]_q!}.
\ee
Since $\pi(\bfz)$ is known the remaining task is to construct $D(\bfz,\bfeta)$
as stated in the theorem by proving \eref{dualityalter}.

\subsubsection{Technical lemmas}

We prove the following lemmas:

\begin{lmm}
\label{LemmaNMbfz}
Consider coordinate sets $\bfx' = \bfx \cup \bfr$ and $\bfy' = \bfy \cup \bfs$.
For $k \notin \bfr$ and $l \notin \bfs$ one has
\bea
\label{f7}
N_{k}(\{\bfx\cup\bfr,\cdot\}) & = & N_{k}(\{\bfx,\cdot\}) + N_{k}(\{\bfr,\cdot\}) \\
\label{f7a}
& = & N_{k}(\{\bfx,\cdot\}) + N(\{\bfr,\cdot\}) - \sum_{i=1}^{N(\{\bfr,\cdot\})} \Theta(k,r_i)\\
\label{f8}
M_{l}(\{\cdot,\bfy\cup\bfs\}) & = & M_{l}(\{\cdot,\bfy\}) + M_{l}(\{\cdot,\bfs\}) \\
\label{f8a}
& = & M_{l}(\{\cdot,\bfy\}) + 
M(\{\cdot,\bfs\}) - \sum_{i=1}^{M(\{\cdot,\bfs\})} \Theta(l,s_i).
\eea
\end{lmm}

{\it Proof:} The function $\Theta(r,x)$ defined in \eref{Thetaxz} satisfies
\bea
\label{f1}
& & \Theta(r,x) = 1 - \Theta(x,r) - \delta_{r,x} \\
\label{f2}
& & \sum_{k=-L+1}^{x-1} \delta_{r,k} = \Theta(r,x), \quad \sum_{k=x+1}^{L} \delta_{r,k} = \Theta(x,r)
\eea

From \eref{partnum}, \eref{occupos} and \eref{NyMx} we have
\bea
\label{f3}
\sum_{i=1}^{N(\bfz)} \sum_{l=r+1}^{L} \delta_{x_i,l} & = &
N(\bfz) - N_{r}(\bfz) - \sum_{i=1}^{N(\bfz)} \delta_{x_i,r} \\
\label{f4}
\sum_{i=1}^{M(\bfz)} \sum_{l=r+1}^{L} \delta_{y_i,l} & = &
M(\bfz) - M_{r}(\bfz) - \sum_{i=1}^{M(\bfz)} \delta_{y_i,r}.
\eea
Specifically, for $N(\bfz)=1$ or $M(\bfz)=1$ resp. we obtain from \eref{f1} - \eref{f4}
\bel{f5}
N_r(\{x,\cdot\}) = \Theta(x,r), \quad M_r(y) = \Theta(y,r),
\ee
and more generally one finds for $\bfz = \{\bfx,\bfy\}$ from \eref{NyMx} and \eref{f2}
\be
\label{f6}
N_{r}(\{\bfx,\cdot\}) = \sum_{i=1}^{N(\bfz)} \Theta(x_i,r) , \quad
M_{r}(\{\cdot,\bfy\}) = \sum_{i=1}^{M(\bfz)} \Theta(y_i,r).
\ee
With \eref{f5} and \eref{f6} one then finds
\be
\label{f6a}
N_{r}(\{\bfx,\cdot\}) = \sum_{i=1}^{N(\bfz)} N_{r}(\{x_i,\cdot\}), \quad
M_{r}(\{\cdot,\bfy\}) = \sum_{i=1}^{M(\bfz)} M_{r}(\{\cdot,y_i\}).
\ee

The first equality \eref{f7} in the lemma then follows from \eref{f6a}. The second equality
\eref{f7a} arises from \eref{f1} and \eref{f6}, bearing in mind that by assumption $k \notin \bfr$. 
The proof of \eref{f8} and \eref{f8a} is analogous. \qed

In particular, for $\bfx = \bfy = \emptyset$ one obtains from Lemma \eref{LemmaNMbfz} 
for $k \notin \bfr$ and $l \notin \bfs$ the inversion formulas
\bea
\label{f9}
N_{k}(\{\bfr,\cdot\}) & = & N(\{\bfr,\cdot\}) - \sum_{i=1}^{N(\{\bfr,\cdot\})} N_{r_i}(\{k,\cdot\})\\
\label{f10}
M_{l}(\{\cdot,\bfs\}) & = & M(\{\cdot,\bfs\}) - \sum_{i=1}^{M(\{\cdot,\bfs\})} M_{s_i}(\{\cdot,l\}).
\eea

We also derive the projector lemma.

\begin{lmm}
\label{projlem}
The tensor occupation operators $\hat{a}_k$, $\hat{b}_k$ act as projectors
\bea
\label{eigenhata}
&& \hat{a}_k \ket{\bfeta} = a_k \ket{\bfeta} = \sum_{i=1}^{N(\bfeta)} \delta_{x_i,k} \ket{\bfeta}\\
\label{eigenhatb}
&& \hat{b}_k \ket{\bfeta} = b_k \ket{\bfeta} = \sum_{i=1}^{M(\bfeta)} \delta_{y_i,k} \ket{\bfeta}
\eea
with the occupation variables $a_k$ and $b_k$ \eref{lonv} (or particle coordinates $x_i$ and $y_i$ respectively)
understood as functions of $\bfeta$ or $\bfz=\bfeta$.
\end{lmm}

{\it Proof:} The first equality in each equation is inherited from the definition of the
projectors
\eref{projection} by multilinearity of the tensor product, the
second equality follows from \eref{occupos}. \qed

Finally we note two combinatorial identities for sums 
over the permutation group $S_n$. One has
\bel{qfacn}
\sum_{\sigma \in S_n} q^{-2 \sum_{j=1}^n \sum_{i=1}^{j-1} \sigma(\Theta(r_i,r_j)) + n(n-1)/2}
= [n]_q! \ q^{-2 \sum_{j=1}^n \sum_{i=1}^{j-1} \Theta(r_j,r_i)},
\ee
which can be proved by induction using $[n]_q = \sum_{k=0}^{n-1} q^{2k-n+1}$,
and 
\bel{permutsum}
\sum_{r_1=-L+1}^L  \sum_{r_n=-L+1}^L f(r_1,\dots,r_n) = \sum_{\vec{r}_n} \sum_{\sigma \in S_n}
f(\sigma(r_1,\dots,r_n))
\ee
for functions that vanish whenever $r_i=r_j$.
Here the sum over $\vec{r}_n$ denotes the summation over the 
Weyl alcove $W^{2L}_n$. 

\subsubsection{Main steps}

After these preparations we go on to prove for $\bfz=\{\bfx,\bfy\}$ the property
\bel{intermediate}
\bra{\bfx,\bfy} S = 
\bra{s} \left( \prod_{i=1}^{N(\{\bfx,\bfy\})}\hat{Q}_{x_i}^A 
\prod_{i=1}^{M(\{\bfx,\bfy\})} \hat{Q}_{y_i}^B \right).
\ee
The matrices on the r.h.s. are the operator form of the functions
$Q_{A,B}^{\cdot}$ defined in \eref{QAB}. According to \eref{dualityalter} proving
\eref{intermediate} proves the theorem.

From the representation \eref{X1mtrafo}, from \eref{eigenhata} and from
the definition \eref{ABxyr}
one finds 
\be 
\bra{\bfx,\bfy} Y_1^-(r) = q^{-A_r(\{\bfx,\bfy\})} \bra{\bfx\cup r,\bfy}.
\ee
By iteration
\be 
\bra{\bfx,\bfy} Y_1^-(r_1) \dots Y_1^-(r_n) = 
q^{-\sum_{j=1}^n A_{r_j}(\{\bfx\cup \bfr_{j-1},\bfy\})} \bra{\bfx\cup \bfr_n,\bfy}
\ee
with the definition $\bfr_0 = \emptyset$.
Since from \eref{ABxyr}, \eref{f7}, \eref{f1} one can write
\be 
A_{r_j}(\{\bfx\cup \bfr_{j-1},\cdot\}) = 2 N_{r_j}(\{\bfx,\cdot\}) + 
2 \sum_{i=1}^{j-1} \Theta(r_i,r_j) - (N(\bfx,\cdot\})+j-1)
\ee
one has
\be 
\sum_{j=1}^n A_{r_j}(\{\bfx\cup \bfr_{j-1},\cdot\}) = 
2 \sum_{j=1}^n N_{r_j}(\{\bfx,\cdot\}) - n N(\{\bfx,\cdot\}) + 
2 \sum_{j=1}^n \sum_{i=1}^{j-1} \Theta(r_i,r_j) - \half n(n-1).
\ee

Now we observe that the term $\sum_{j=1}^n (2 N_{r_j}(\{\bfx,\cdot\}) - N(\{\bfx,\cdot\}))$ 
is invariant under permutations of the coordinates $r_j$.
Using the fact that $(Y_1^-(r))^2 = 0$ and the combinatorial properties \eref{qfacn} and \eref{permutsum} 
then yields
\be 
\bra{\bfx,\bfy} \frac{(Y_1^-)^n}{[n]_q!} = \sum_{\vec{r}_n}
q^{-\sum_{j=1}^n (2 N_{r_j}(\{\bfx,\bfy\}) - N(\{\bfx,\bfy\}))} \bra{\bfx\cup \bfr_n,\bfy}.
\ee
The next step is to invoke Lemma \eref{LemmaNMbfz} to express $N_{r_j}(\{\bfx,\cdot\})$ 
in terms of single-particle step functions $N_{x_i}(\{r_j,\cdot\})$ with inverted arguments.
This initiates the following chain of equalities for the exponent 
$E_{\bfr}(\{\bfx,\cdot\}) := - \sum_{j=1}^n[ 2N_{r_j}(\{\bfx,\cdot\}) - N(\{\bfx,\cdot\})]$
of $q$:
\bea
E_{\bfr}(\{\bfx,\cdot\}) & = & - \sum_{j=1}^n\left( N(\{\bfx,\cdot\}) - 
2 \sum_{i=1}^{N(\{\bfx,\cdot\})} N_{x_i}(\{r_j,\cdot\})\right) \nonumber \\
& = & \sum_{j=1}^n \sum_{i=1}^{N(\{\bfx,\cdot\})} ( 2  N_{x_i}(\{r_j,\cdot\}) - 1)  \\
& = & \sum_{i=1}^{N(\{\bfx,\cdot\})} A_{x_i}(\{\bfr,\cdot\}) . \nonumber
\eea

Since trivially
\be 
\bra{\bfx\cup \bfr_n,\bfy} = \bra{\bfx\cup \bfr_n,\bfy} a_{x_1} \dots a_{x_N}
\ee
one arrives at
\bea
\bra{\bfx,\bfy} \frac{(Y_1^-)}{[n]_q!} & = & \sum_{\vec{r}_n} 
\left( \prod_{i=1}^{N(\{\bfx,\bfy\})} q^{A_{x_i}(\{\bfr,\bfy\})} a_{x_i} \right) 
\bra{\bfx\cup \bfr_n,\bfy} \nonumber \\
& = & \sum_{\vec{r}_{\tilde{n}}} 
\left( \prod_{i=1}^{N(\{\bfx,\bfy\})} q^{A_{x_i}(\{\bfr,\bfy\})} a_{x_i} \right) \bra{\bfr_{\tilde{n}},\bfy}. 
\eea
where the summation in the second equality has been changed to the
extended Weyl alcove with $\tilde{n} = N(\{\bfx,\bfy\})+n$. This is possible since
the product of indicators $a_{x_1} \dots a_{x_N}$ cancels all terms not belonging to the original
Weyl alcove $n$.

Next one uses the definitions \eref{ABxyr}, \eref{QAB} and 
the projector property \eref{eigenhata} to express $A_{x_i}(\{\bfr,\bfy\})$ in terms
of projectors. The result is
\bea
\bra{\bfx,\bfy} \frac{(Y_1^-)^n}{[n]_q!} & = & \sum_{\vec{r}_{\tilde{n}}} 
\bra{\bfr_{\tilde{n}},\bfy} \left( \prod_{i=1}^{N(\{\bfx,\bfy\})} 
 q^{\sum_{j=1}^{x_i-1} \hat{a}_j - \sum_{j=x_i+1}^{L} \hat{a}_j} \hat{a}_{x_i} \right) \nonumber \\
& = & \sum_{\vec{r}_{\tilde{n}}} \bra{\bfr_{\tilde{n}},\bfy} \prod_{i=1}^{N(\{\bfx,\bfy\})}  \hat{Q}^A_{x_i}.
\eea

Using the commutator relation $[Y_1^-,Y_2^+]=0$ 
\eref{Uqglncomm3} and going through similar steps yields
\be
\bra{\bfx,\bfy} \frac{(Y_1^-)^n}{[n]_q!} \frac{(Y_2^+)^m}{[m]_q!} =
\sum_{\vec{r}_{\tilde{n}}} \sum_{\vec{s}_{\tilde{m}}} \bra{\bfr_{\tilde{n}},\bfs_{\tilde{m}}}
\left( \prod_{i=1}^{N(\{\bfx,\bfy\})} \hat{Q}^A_{x_i}\prod_{i=1}^{M(\{\bfx,\bfy\})} \hat{Q}^B_{y_i} \right)
\ee
with $\tilde{m} = M(\{\bfx,\bfy\})+m$ and
\be
\hat{Q}^B_y = q^{-\sum_{k=-L+1}^{y-1} \hat{b}_k + \sum_{k=y+1}^{L} \hat{b}_k} \hat{b}_y.
\ee

Observe now that the summation on the r.h.s involves only the vector 
$\bra{\bfr_{\tilde{n}},\bfs_{\tilde{m}}}$.
Since the summation is over the Weyl alcove one has
\be 
\sum_{\vec{r}_{\tilde{n}}} \sum_{\vec{s}_{\tilde{m}}} 
\bra{\bfr_{\tilde{n}},\bfs_{\tilde{m}}} = \bra{s_{\tilde{n},\tilde{m}}}.
\ee
Therefore
\be
\bra{\bfx,\bfy} \frac{(Y_1^-)^n}{[n]_q!} \frac{(Y_2^+)^m}{[m]_q!} = \bra{s_{\tilde{n},\tilde{m}}} 
\left( \prod_{i=1}^{N(\{\bfx,\bfy\})} \hat{Q}^A_{x_i} \prod_{i=1}^{M(\{\bfx,\bfy\})} \hat{Q}^B_{y_i} \right).
\ee
Notice that on the r.h.s. the only dependence on the $n$ and $m$ is in the summation vector
\bra{s_{\tilde{n},\tilde{m}}} for the sector with $\tilde{n}=N(\{\bfx,\bfy\})+n$ particles of type $A$ and
$\tilde{m}=M(\{\bfx,\bfy\})+m$ particles of type $B$.

The final step is to take the double sum \eref{S}. Terms such that 
$n+m > 2L - N(\{\bfx,\bfy\}) - M(\{\bfx,\bfy\})$
are zero since $\bra{s_{N,2L-N}} Y_1^- = \bra{s_{N,2L-N}} Y_2^+ = 0$, corresponding to the exclusion
principle that forbids creating configurations with more than $2L$ particles on $\Lambda$. 
On the other hand,  
\be 
\bra{s_{n,m}} \left(\prod_{i=1}^{N(\{\bfx,\bfy\})}  \hat{Q}^A_{x_i} 
\prod_{i=1}^{M(\{\bfx,\bfy\})}  \hat{Q}^B_{y_i} \right) = 0 
\mbox{ if } n<N(\{\bfx,\bfy\}) \mbox{ or } m<M(\{\bfx,\bfy\}),
\ee
due to the projectors contained in the operators $\hat{Q}^{A,B}_{.}$
This yields
\eref{intermediate} with
\bel{dualitymatrixform}
\hat{Q}_\bfz = \hat{Q}^A_\bfx \hat{Q}^B_\bfy 
\ee
and proves Theorem \eref{Theo:Duality} 
by taking the scalar product with $\bra{s}$ and $\ket{\bfeta}$. \qed

\subsection{Proof of Theorem \eref{Theo:Sumrule}}

The first equality follows from duality and ergodicity by taking the limit 
$t \to \infty$ in the expectation \eref{coroQ} with coordinate sets $\bfz'$ 
representing configurations with
$N'$ particles of type $A$ and $M'$ particles of type $B$.

We also have from \eref{coroQ} by taking $P_0=\pi^\ast_{N,M}$, i.e., by considering
the canonical equilibrium distribution for $N$ particles of type $A$ and 
$M$ particles of type $B$,
\be
\exval{ Q_{\bfz}}_{\pi^\ast_{N,M}} =  \bra{\bfz} \rme^{-Ht} \ket{V_{N,M}}.
\ee
with a vector 
\be 
\ket{V_{N,M}} = \sum_{\bfz'} \exval{ Q_{\bfz'}}_{\pi^\ast_{N,M}} \ket{\bfz'}.
\ee
that has support in the subspace corresponding to configurations with
$N'$ particles of type $A$ and $M'$ particles of type $B$.

Because of stationarity the l.h.s. does not depend on time. This implies by ergodicity that
$\ket{V_{N,M}}$ on the r.h.s. must be proportional to the (unique) stationary probability vector
for configurations in $\S^{2L}_{N',M'}$. Hence
\bel{Qlambda}
\exval{ Q_{\bfz}}_{\pi^\ast_{N,M}} = \lambda_{N,M}^{N',M'} \pi^\ast_{N',M'} (\bfz) 
\quad \forall \bfz \in \mathbb{S}^{2L}_{N',M'}
\ee
where $\lambda_{N,M}^{N',M'}$ is some constant with $N' = N(\bfz) = N(\bfz')$ and $M' = M(\bfz)= M(\bfz')$.
This proves the second equality in the theorem. \qed

\section*{Acknowledgements}
This work was supported by DFG and by CNPq through the grant 307347/2013-3.
GMS thanks the University of S\~ao Paulo,
where part of this work was done, for kind hospitality.

\end{document}